\documentclass[]{article}

\usepackage{amssymb, amsbsy, amsmath, amsfonts, amssymb, amscd, mathrsfs}

\usepackage[english]{babel}
\usepackage[T1]{fontenc}
\usepackage{indentfirst}

\usepackage{ifpdf}
\ifpdf
\usepackage[pdftex]{color, graphicx}
\else
\usepackage[dvips]{color, graphicx} 
\fi

\usepackage{caption}

\newtheorem{statement}{}[section]
\newtheorem{theorem}[statement]{Theorem}

\newtheorem{proposition}[statement]{Proposition}

\newtheorem{corollary}[statement]{Corollary}

\newcommand\C{\mathbb C}

\newcommand\T{\mathbb T}
\newcommand\D{\mathbb D}

\newcommand\e{{\rm e}}

\newcommand\eps{\varepsilon}
\newcommand\ind{{\rm 1\kern-.30em I}}
\newcommand\qed{\hfill $\square$}
\renewcommand \Re{{\mathfrak R}{\rm e}\,}

\let\phi=\varphi

\title{Compact composition operators on Hardy-Orlicz and Bergman-Orlicz spaces
\footnote{
These results come from joint works with P. Lef\`evre, H. Queff\'elec and L. Rodr{\'\i}guez-Piazza 
(\cite{CompOrli}, \cite{LLQR-N}, \cite{LLQR-B}). It is an expanded version of the conference I gave at the ICM 
satellite conference \emph{Functional Analysis and Operator Theory}, held in Bangalore, India, 8--11 august 2010.} 
}

\author{Daniel Li}
\date{\footnotesize \today}

\begin{document}

\maketitle

\noindent{\bf Abstract.} 
It is known, from results of B. MacCluer and J. Shapiro (1986), that every composition 
operator which is compact on the Hardy space $H^p$, $1 \leq p < \infty$, is also compact on the Bergman space 
${\mathfrak B}^p = L^p_a (\D)$. In this survey, after having described the above known results, we consider 
Hardy-Orlicz $H^\Psi$ and Bergman-Orlicz ${\mathfrak B}^\Psi$ spaces, characterize the compactness of their composition 
operators, and show that there exist Orlicz functions for which there are composition operators which are compact on 
$H^\Psi$ but not on ${\mathfrak B}^\Psi$.\par
\medskip

\medskip


\section{Introduction}

Let $\D =\{z \in \C\,;\ |z| < 1\}$ be the open unit disk of the complex plane. For $1\leq p < \infty$, consider the 
Hardy space 
\begin{displaymath}
H^p = \{ f \colon \D \to \C \,;\ f \text{ analytic and } \|f \|_{H^p} < +\infty\}\,,
\end{displaymath}
where
\begin{displaymath}
\|f \|_{H^p} = \sup_{r < 1} \bigg[\frac{1}{2\pi} \int_0^{2\pi} |f (r\e^{it})|^p\, dt\bigg]^{1/p},
\end{displaymath}
and the Bergman space ${\mathfrak B}^p$ (otherwise denoted by $A^p$ or $L^p_a$) 
\begin{displaymath}
{\mathfrak B}^p = \{f \colon \D \to \C\,; \ f \text{ analytic and } f \in L^p (\D, {\cal A})\}\,,
\end{displaymath}
whose norm is defined by
\begin{displaymath}
\| f \|_{{\mathfrak B}^p} = \bigg[\int_\D |f (z)|^p\, d{\cal A} (z) \bigg]^{1/p}\,,
\end{displaymath}
where $d{\cal A} = \frac{dx dy}{\pi}$ is the normalized area measure on $\D$.\par
\bigskip

Every analytic self-map $\phi \colon \D \to \D$ (such a function is also known as a Schur function, or function of 
the Schur-Agler class) defines a \emph{composition operator} $f \mapsto C_\phi (f) = f \circ \phi$ which is a 
\emph{bounded} linear operator 
\begin{displaymath}
C_\phi \colon H^p \to H^p \qquad\quad  \text{resp.} \quad C_\phi \colon {\mathfrak B}^p \to {\mathfrak B}^p,
\end{displaymath}
thanks to \emph{Littlewood's subordination principle} (see \cite{Duren}, Theorem~1.7). The function $\phi$ is called the 
\emph{symbol} of $C_\phi$.

\subsection{Compactness}

The compactness of composition operators on Bergman spaces had been characterized in 1975 by D. M. Boyd (\cite{Boyd}) 
for $p = 2$ and by B. MacCluer and J. Shapiro in 1986 for the other $p$'s (\cite{McCluer-Shapiro}, Theorem~3.5 and 
Theorem~5.3):
\begin{theorem} [Boyd (1976), MacCluer-Shapiro (1986)]\label{Theo McCluer-Shapiro}
Let $\phi$ be an analytic self-map of $\D$. Then, for $1 \leq p < \infty$, one has:
\begin{equation}\label{Bergman angulaire}
C_\phi \colon {\mathfrak B}^p \to {\mathfrak B}^p \quad \text{compact} \qquad \Longleftrightarrow \qquad 
\lim_{|z| \to 1} \frac{1 - |\phi (z)|}{1 - |z|} = + \infty\,.
\end{equation}
\end{theorem}
That means that $\phi (z)$ approaches the boundary of $\D$ less quickly that $z$.
That means also (but we do not need this remark in the sequel) that $\phi$ has no finite angular derivative on 
$\partial \D$ ($\phi$ has a finite angular derivative at $\omega \in \partial \D$ if the angular limits 
$\phi^\ast (\omega) = \lim_{z \to \omega} \phi (z)$ and 
$\phi' (\omega) = \angle \lim_{z \to \omega} \frac{\phi (z) - \phi^{\ast} (\omega)}{z - \omega}$ exist, with $|\phi^\ast (\omega)| = 1$; by the 
Julia-Carath\'eodory Theorem, the non-existence of such a derivative is equivalent to the right-hand side of 
\eqref{Bergman angulaire}: see \cite{Shapiro}, \S~4.2). We actually prefer to write \eqref{Bergman angulaire} in 
the following way:
\begin{equation}\label{Bergman angulaire bis}
C_\phi \colon {\mathfrak B}^p \to {\mathfrak B}^p \quad \text{compact} \qquad \Longleftrightarrow \qquad 
\lim_{|z| \to 1} \frac{1 - |z|}{1 - |\phi (z)|} = 0\,.
\end{equation}\par
\medskip

On the other hand, it is not difficult to see (\cite{Shapiro-Taylor}, Theorem~2.1, or \cite{Shapiro}, \S~3.5):
\begin{proposition}[Shapiro-Taylor (1973)]\label{Shapiro-Taylor}
Let $\phi$ be an analytic self-map of $\D$. Then, for $1 \leq p < \infty$, one has:
\begin{equation}\label{theo Shapiro-Taylor}
C_\phi \colon H^p \to H^p \quad \text{compact} \qquad \Longrightarrow \qquad 
\lim_{|z| \to 1} \frac{1 - |z|}{1 - |\phi (z)|} = 0\,.
\end{equation}
\end{proposition}
(this is actually an equivalence when $\phi$ is univalent (\cite{Shapiro}, \S~3.2), or more generally boundedly-valent, 
but there are Blaschke products for which the converse of \eqref{theo Shapiro-Taylor} is not true: see \cite{Shapiro}, 
\S~10.2, or \cite{LLQR-Rev}, Theorem~3.1, for a more general result, with a simpler proof).\par\medskip

Hence:
\begin{corollary}\label{comparaison cas classique}
For $1 \leq p < \infty$, one has:
\begin{displaymath}
C_\phi \colon H^p \to H^p \quad \text{compact} \qquad \Longrightarrow \qquad 
C_\phi \colon {\mathfrak B}^p \to {\mathfrak B}^p \quad \text{compact}. 
\end{displaymath}
The converse is not true.
\end{corollary}

\subsection{Goal}

Our goal is to replace the classical Hardy spaces $H^p$ and Bergman spaces ${\mathfrak B}^p$ by Hardy-Orlicz spaces 
$H^\Psi$ and Bergman-Orlicz spaces ${\mathfrak B}^\Psi$ and compare the compactness of composition operators 
$C_\phi \colon H^\Psi \to H^\Psi$ and $C_\phi \colon {\mathfrak B}^\Psi \to {\mathfrak B}^\Psi$. We shall detail that 
in Section~\ref{section comparaison}.

\subsection{Some comments}\label{some comments}

To prove Proposition~\ref{Shapiro-Taylor}, J. H. Shapiro and P. D. Taylor used the following result 
(\cite{Shapiro-Taylor}, Theorem~6.1):

\begin{theorem}[Shapiro-Taylor (1973)]\label{independance de p}
If the composition operator $C_\phi$ is compact on $H^{p_0}$ for some $1 \leq p_0 < \infty$, then it is compact on 
$H^p$ for all $1 \leq p < \infty$.
\end{theorem}

\noindent{\bf Proof.} First, by Montel's Theorem, $C_\phi$ is compact on $H^p$ if and only if 
$\|C_\phi (f_n)\|_p$ converges to $0$ for every sequence $(f_n)$ in the unit ball of $H^p$ which converges uniformly 
to $0$ on compact subsets of $\D$. One then uses Riesz's factorization Theorem: if $f_n$ is in  
$H^p$, we can write $f_n = B_n g_n$, where $B_n$ is a Blaschke product and $g_n$ has no zero in $\D$. By Montel's 
Theorem, we may assume that $(g_n)_n$ converge uniformly on compact sets of $\D$. Setting 
$h_n = g_n^{p/p_0}$, we get a function which is in the unit ball of $H^{p_0}$. Since $(g_n)_n$ converges uniformly on 
compact sets of $\D$, so does $(h_n)_n$. Its limit $h$ belongs also to the unit ball of $H^{p_0}$ and, 
by the compactness of $C_\phi$, $\big( C_\phi (h_n) \big)_n$ converges to $C_\phi (h)$ in $H^{p_0}$.\par  
Now the compactness of $C_\phi$ on $H^{p_0}$ implies that $|\phi^\ast| < 1$ almost everywhere ($\phi^\ast$ is the 
boundary values function of $\phi$). In fact, let 
$P_n (z) = z^n$, since $(P_n)$ is in the unit ball of $H^{p_0}$ and converges uniformly to $0$ on compact subsets of 
$\D$, one has $\| C_\phi (P_n)\|_{p_0} \to 0$. But, if $E_\phi =\{\xi \in \partial \D\,;\ |\phi^\ast (\xi)| = 1\}$, 
one has $\|C_\phi (P_n)\|_{p_0}^{p_0} \geq \int_{E_\phi} |\phi^\ast (\xi)|^{n p_0}\,dm (\xi) \geq m (E_\phi)$, where 
$m$ is the normalized Lebesgue measure on $\partial \D$. Hence $m (E_\phi) = 0$.\par
It follows that the sequence $\big((h_n \circ \phi)^\ast\big)_n = (h_n \circ \phi^\ast)_n$ converges almost everywhere 
to $h \circ \phi^\ast$ on $\partial \D$. Since $\big( C_\phi (h_n) \big)_n$ converges in the norm of $L^{p_0} (\T)$, Vitali's convergence 
Theorem gives:
\begin{displaymath}
\lim_{m (E) \to 0} \sup_n \int_E |h_n \circ \phi^\ast|^{p_0}\,dm = 0\,.
\end{displaymath}
But
\begin{displaymath}
\int_E |f_n \circ \phi^\ast|^p\,dm \leq \int_E |g_n \circ \phi^\ast|^p\,dm = 
\int_E |h_n \circ \phi^\ast|^{p_0}\,dm\,,
\end{displaymath}
so Vitali's convergence Theorem again gives $\| f_n \circ \phi \|_p \to 0$, since $f_n \circ \phi^\ast$ tends to $0$ {\it a.~e.} on $\partial \D$.\qed
\bigskip

Actually, the proof of Proposition~\ref{Shapiro-Taylor} can be made without using Theorem~\ref{independance de p}, but 
we gave it to see that Riesz's factorization Theorem is the main tool.
\par\medskip 

\noindent{\bf Proof of Proposition~\ref{Shapiro-Taylor}.} For every $z \in \D$, the evaluation map 
$e_z \colon f \in H^p \mapsto f (z)$ is a continuous linear form and $\| e_z \| \leq \frac{2^{1/p}}{(1 - |z|)^{1/p}}$ 
(see \cite{Duren}, lemma in \S~3.2, page 36). But actually $\|e_z\| = \frac{1}{(1 - |z|^2)^{1/p}}\,$. 
Indeed, let $u_z (\zeta) = \big(\frac{1 - |z|}{1 - \bar{z} \zeta}\big)^{2/p}$, $|\zeta| < 1$. Then, using the 
Parseval formula, $\|u_z \|_p^p = \frac{1 - |z|}{1 + |z|}\,$. Therefore
\begin{displaymath}
\|e_z\| \geq \frac{|u_z (z)|}{\|u_z\|_p} 
\geq \frac{\frac{1}{\hskip 3pt (1 + |z|)^{2/p}}}{\hskip 3pt \big(\frac{1 - |z|}{1 + |z|}\big)^{1/p}} 
= \frac{1}{(1 - |z|^2)^{1/p}}\,\cdot
\end{displaymath}
On the other hand, it is clear, by the Cauchy-Schwarz inequality, that $|h (z)| \leq 1/ (1 - |z|^2)^{1/2}$ for every 
$h$ in the unit ball of $H^2$; hence, if $f$ is in the unit ball of $H^p$, and we write $f = Bg$, where $B$ is the 
Blaschke product associated to the zeroes of $f$, we get, since $h = g^{p/2}$ is in the unit ball of $H^2$: 
$|f(z)| \leq |g (z)| = |h(z)|^{2/p} \leq 1/ (1 - |z|^2)^{1/p}$. 
\par\medskip

\noindent
Now, $\frac{e_z}{\|e_z\|} \mathop{\longrightarrow}\limits^{w^\ast}_{|z| \to 1} 0$ in $(H^p)^\ast$, because 
$e_z (P)(1 - |z|)^{1/p} = P (z) (1 - |z|)^{1/p} \mathop{\longrightarrow}\limits_{|z| \to 1} 0$ for every polynomial 
$P$. Since $C_\phi$ is compact, its adjoint is also compact; hence 
$\|C_\phi^\ast (e_z/\|e_z\|) \| \mathop{\longrightarrow}\limits_{|z| \to 1} 0$, and that gives the result since 
$\|C_\phi^\ast (e_z/\|e_z\|) \| = \frac{\| e_{\phi (z)}\|}{\|e_z\|}  = \Big( \frac{1 - |z|^2}{1 - |\phi (z)|^2}\Big)^{1/p} 
\geq \frac{1}{2^{1/p}} \Big(\frac{1 - |z|}{1 - |\phi (z)|} \Big)^{1/p}$.\qed
\par\bigskip

For Bergman spaces, the necessary condition of compactness in Theorem~\ref{Theo McCluer-Shapiro} follows the same 
lines as in the Hardy case. B. D. MacCluer and J. H. Shapiro (\cite{McCluer-Shapiro}, Theorem~5.3) proved the 
sufficient condition, in showing that the compactness of $C_\phi$ on one of the Bergman spaces ${\mathfrak B}^{p_0}$, 
$1 \leq p_0 < \infty$, implies its compactness for all the Bergman spaces ${\mathfrak B}^p$, $1 \leq p < \infty$ 
(see Theorem~\ref{compactness classical} below), and then used Boyd's result for $p = 2$ (actually, they gave a new 
proof of Boyd's result). To do that, since there is no Bergman version of Riesz's factorization Theorem, they had to 
use another tool and they used the notion of \emph{Carleson measure}, that we shall develop in 
the next section. Before that, let us give a  proof of Boyd's result. We follow \cite{Zhu}, Theorem~10.3.5.\par
\medskip

\noindent{\bf Proof of Theorem~\ref{Theo McCluer-Shapiro}} (for $p = 2$). We only have to show that  
$\lim\limits_{|z| \to 1} \frac{1 - |z|}{1 - |\phi (z)|} = 0$ implies the compactness of $C_\phi$ on 
${\mathfrak B}^2$.\par 
We may assume that $\phi (0) = 0$. Indeed, if $\phi (0) = a$, let $\psi = \phi_a \circ \phi$, with 
$\phi_a (z) = \frac{a - z}{1 - \bar{a} z}$\,. One has $\psi (0) = 0$ and $C_\phi = C_\psi \circ C_{\phi_a}$, since $\phi = \phi_a \circ \psi$. 
Hence the compactness of $C_\psi$ implies the one of $C_\phi$. Moreover, one has $|\psi (z)| \leq |\phi (z)|$, so that the condition 
$\lim_{|z| \to 1} \frac{1 - |z|}{1 - |\phi (z)|} = 0$ implies that $\lim_{|z| \to 1} \frac{1 - |z|}{1 - |\psi (z)|} = 0$.
\par
Let $(f_n)$ be a sequence in the unit ball of ${\mathfrak B}^2$ which converges to $0$ uniformly on compact subsets of 
$\D$. Then so does $(f'_n)$. Using Taylor expansion, one has a 
constant $ C > 0$ such that:
\begin{displaymath}
\int_\D |f (z)|^2\, d{\cal A} (z) \leq C\,\bigg[ |f( 0)|^2 + \int_\D (1 - |z|^2)^2 |f' (z)|^2\,d{\cal A}(z)\bigg]
\end{displaymath}
for every analytic function $f \colon \D \to \C$. It follows that
\begin{displaymath}
\|C_\phi (f_n)\|_{{\mathfrak B}^2}^2 \leq C\, \bigg[ |(f_n \circ \phi) (0)|^2 
+ \int_\D (1 - |z|^2)^2 |(f_n \circ \phi)' (z)|^2\,d{\cal A}(z)\bigg]
\end{displaymath}
for every $n \geq 1$. Since $f_n [\phi (0)] \mathop{\longrightarrow}\limits_{n \to \infty} 0$, it remains to show 
that the integral tends to $0$.\par
For every $\eps > 0$, we may take, by hypothesis, some $\delta > 0$, with $\delta < 1$, such that 
$\frac{1 - |z|^2}{1 - |\phi (z)|^2} \leq \eps$ for $\delta \leq |z| < 1$. This implies that:
\begin{align*}
\int_\D (1 - |z|^2)^2  |(f_n & \circ \phi)' (z)|^2  \,d{\cal A}(z) \\
 \leq  \int_{|z| < \delta} & \hskip - 5pt (1 - |z|^2)^2 |(f'_n \circ \phi) (z)|^2 |\phi' (z)|^2\,d{\cal A}(z) \\ 
& + \eps  \int_{\delta \leq |z| < 1} \hskip -8 pt
(1 - |\phi (z)|^2) (1 - |z|^2) |\phi' (z)|^2 |(f'_n \circ \phi) (z)|^2 \,d{\cal A}(z)\,.
\end{align*}
Denote by $I_n$ the first integral and by $J_n$ the second one. Since $f'_n [\phi (z)]$ tends to $0$ uniformly for 
$|z| \leq \delta$, $I_n$ tends to $0$. It remains to show that the sequence $(J_n)$ is bounded. Since 
$1 - |z|^2 \leq 2 \log 1/|z|$, the change of variable formula (see \cite{Shapiro}, p.~179, or \cite{Zhu}, 
Proposition~10.2.5) gives:
\begin{align*}
J_n 
& \leq 2 \int_\D (1 - |\phi (z)|^2) |\phi' (z)|^2 |(f'_n \circ \phi) (z)|^2 \log \frac{1}{|z|}\, d{\cal A} (z) \\
& = 2 \int_\D (1 - |w|^2) |f'_n (w)|^2 N_\phi (w)\, d{\cal A} (w).
\end{align*}
Since $\phi (0) = 0$, Littlewood's inequality (see \cite{Shapiro}, \S~10.4) reads as $N_\phi (w) \leq \log 1/|w|$ for $w \neq 0$. We get, since 
$2 \log 1/|w| \approx 1 - |w|^2$ as $|w| \to 1$, a constant $C > 0$ such that:
\begin{displaymath}
J_n \leq C\,\int_\D (1 - |w|^2)^2 |f'_n (w)|^2\, d{\cal A} (w)\,,
\end{displaymath}
and that ends the proof, since one easily sees (using Taylor expansion) that this last integral is less or equal 
than $\int_\D |f_n (w)|^2\,d{\cal A} (w) \leq 1$.\qed
\goodbreak

\subsection{Carleson measures}

If $\D$ is the open unit disk of the complex plane and $\T = \partial \D$ is the unit circle, we denote by 
${\cal A}$ the normalized area measure on $\D$ and by $m$ the normalized Lebesgue measure on $\T$.\par\smallskip

The \emph{Carleson window $W (\xi, h)$} centered at $\xi \in \T$ and of size $h$, $0 < h \leq 1$, is the set 
\begin{displaymath}
W (\xi, h) = \{ z \in \D\,;\ |z| > 1 - h \quad \text{and } \quad |\arg (z \bar{\xi})| < \pi h\}\,.
\end{displaymath}

\begin{figure}[t]
\begin{center}
\includegraphics[width=0.4\linewidth]{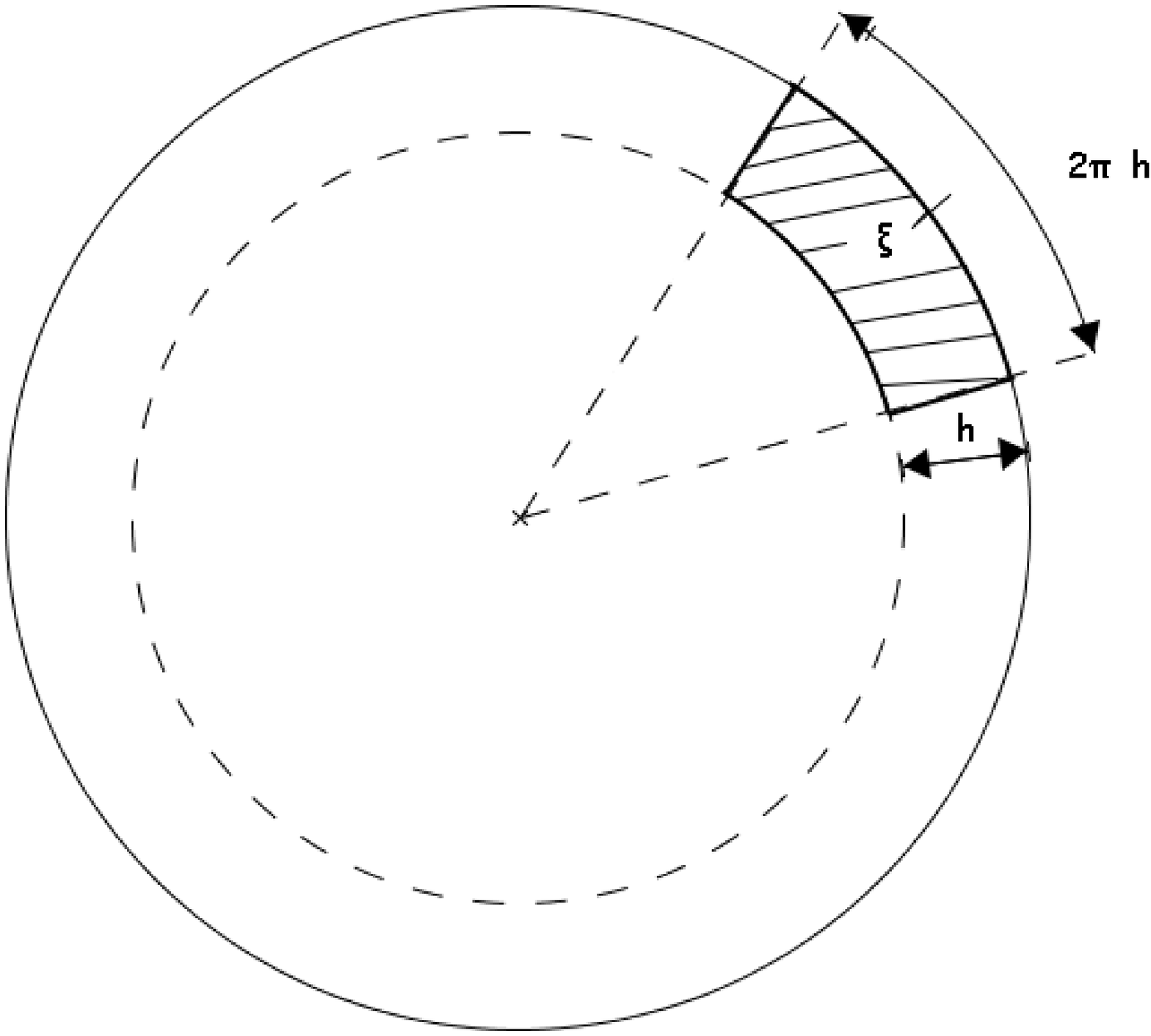}
\end{center}
\end{figure}
The point is that the two dimensions of the window $W (\xi, h)$ are proportional.\par
\smallskip

An \emph{$\alpha$-Carleson measure} (\emph{Carleson measure} if $\alpha = 1$) is a measure $\mu$ on $\D$ such that: 
\begin{displaymath}
\rho_\mu (h) := \sup_{|\xi| = 1} \mu [ W (\xi, h)] = O\,(h^\alpha)\,.
\end{displaymath}
We say that $\rho_\mu$ is the \emph{Carleson function} of $\mu$.\par

We denote by ${\cal A}_\phi$ the pull-back measure of ${\cal A}$ by $\phi$ and by $m_\phi$ the pull-back measure of $m$ by $\phi^\ast$, 
$\phi^\ast$ being the boundary values function of $\phi$. Recall that these pull-back measures are defined by $m_\phi (E) = m[{\phi^\ast}^{- 1} (E)]$ 
and ${\cal A}_\phi (E) = {\cal A} [\phi^{-1} (E)]$ for every Borel set $E$ of $\partial \D$ or of $\D$, respectively. We write 
$\rho_\phi$ and $\rho_{\phi, 2}$ the Carleson functions of $m_\phi$ and ${\cal A}_\phi$ respectively. We call them the \emph{Carleson function of $\phi$} 
and the \emph{Carleson function of order $2$ of $\phi$}.\par
\medskip

Carleson's Theorem (see \cite{Duren}, Theorem~9.3) says that, for $p < \infty$, the inclusion map 
$I_\mu \colon f \in H^p \mapsto f \in L^p (\mu)$ is defined and bounded if and only if $\mu$ is a Carleson measure. 
A Bergman version has been proved by W. W. Hastings in 1975 (\cite{Hastings}): 
$J_\mu \colon f \in {\mathfrak B}^p \mapsto f \in L^p (\mu)$ is defined and bounded if and only if $\mu$ is a 
$2$-Carleson measure. Now, composition operators $C_\phi \colon H^p \to H^p$, resp. 
$C_\phi \colon {\mathfrak B}^p \to {\mathfrak B}^p$, may be seen as inclusion maps 
$I_\phi \colon H^p \to L^p(\D, m_\phi)$, resp. $J_\phi \colon {\mathfrak B}^p \to L^p (\D, {\cal A}_\phi)$, since:
\begin{displaymath}
\|C_\phi (f)\|_{H^p}^p = \int_\D |f|^p\,dm_\phi \qquad  \text{and} \qquad 
\|C_\phi (f)\|_{{\mathfrak B}^p}^p = \int_\D |f|^p\,d{\cal A}_\phi\,.
\end{displaymath} 
Hence the continuity of the composition operator $C_\phi$, both on $H^p$ and 
${\mathfrak B}^p$, implies that $m_\phi$ is a Carleson measure and ${\cal A}_\phi$ is a $2$-Carleson measure. It should 
be stressed that in the Hardy case, $m_\phi$ is a measure on $\overline{\D}$, and not on $\D$ in general, so we have 
to adapt the previous notations in this case.\par
\smallskip

For compactness, one has the following result:
\goodbreak
\begin{theorem}[MacCluer (1985), MacCluer-Shapiro (1986)]\label{compactness classical}\hfill\par
\noindent For $1 \leq p < \infty$, and every analytic self-map $\phi$ of $\D$, one has:
\begin{center}
$C_\phi \colon H^p \to H^p$ \quad compact \qquad $\Longleftrightarrow$ \qquad $\rho_\phi (h) = o\,(h)$, \quad as $h \to 0$,  
\end{center}
and
\begin{center}
\quad $C_\phi \colon {\mathfrak B}^p \to {\mathfrak B}^p$ \quad compact \qquad $\Longleftrightarrow$ \qquad 
$\rho_{\phi, 2} (h) = o\,(h^2)$, \quad as $h \to 0$.
\end{center}
\end{theorem}

Recall that when the composition operator $C_\phi \colon H^p \to H^p$ is compact, one has $|\phi^\ast| < 1$ 
{\it a.e.}, and hence $m_\phi$ is supported by $\D$. Since these characterizations do not depend on $p$, 
one recovers Theorem~\ref{independance de p} and get its Bergman counterpart.\par
\smallskip

We shall see in the next section how this theorem changes when we replace the classical Hardy and Bergman spaces 
by their Orlicz generalizations.\par\medskip

\noindent{\bf Proof of Theorem~\ref{compactness classical}.} We only prove the Bergman case; the Hardy case being 
analogous.\par
1) Assume that $C_\phi$ is compact on ${\mathfrak B}^p$. Consider, for every $a \in \D$, the Berezin kernel 
\begin{displaymath}
H_a  = \frac {(1 - |a|^2)^2} {|1 - \overline{a} z|^4} \,\cdot
\end{displaymath}
One has $\| H_a \|_{{\mathfrak B}^1} = 1$ and 
\begin{displaymath}
\| H_a \|_\infty = \frac {(1 - |a|^2)^2} {(1 - |a|)^4} = \frac {(1 + |a|)^2} {(1 - |a|)^2 } 
\leq \frac{4}{(1 - |a|)^2} \,;
\end{displaymath}
hence, writing $a = (1 - h)\xi$, $0 < h \leq 1$, $|\xi| = 1$, we get $\|H_a\|_{{\mathfrak B}^p} \leq (4/h^2)^{1 - \frac{1}{p}}$ and the function 
$f_a = \big(\frac{h^2}{4}\big)^{1- \frac{1}{p}} H_a$ is in the unit ball of ${\mathfrak B}^p$. Moreover, 
$f_a$ tends to $0$ as $|a| \to 1$ uniformly on compact subsets of $\D$. Since $C_\phi$ is compact on ${\mathfrak B}^p$, 
one has $\|C_\phi (f_a) \|_{{\mathfrak B}^p} \mathop{\longrightarrow}\limits_{|a|\to 1} 0$. 
But it is easy to see that $|1 - \bar{a}z| \leq 5h$ when 
$z \in W (\xi, h)$. Hence $|f_a (z)| \geq C_p/h^{2/p}$ when $z \in W (\xi, h)$ and 
\begin{displaymath}
\|C_\phi (f_a) \|_{{\mathfrak B}^p}^p 
= \int_\D |f_a \circ \phi|^p\,d{\cal A} \geq \int_{W (\xi, h)} |f_a|^p\,d{\cal A}_\phi  
\geq \frac{C_p^p}{h^2}\,{\cal A}_\phi [W (\xi, h)]\,. 
\end{displaymath}
Hence $\rho_{\phi, 2} (h) = o\,(h^2)$.\par\smallskip

2) Conversely, assume that $\rho_{\phi, 2} (h) = o\,(h^2)$, and let $(f_n)$ be a sequence in the unit ball of 
${\mathfrak B}^p$ converging uniformly to $0$ on compact subsets of $\D$, and $\eps > 0$.\par 

For every measure $\mu$ on $\D$, let $K_{\mu, 2} (h) = \sup_{0 < t \leq h} \frac{\rho_\mu (t)}{t^2}$ and 
$K_{\phi, 2} (h) = K_{{\cal A}_\phi, 2} (h)$.\par
By hypothesis, there is a $\delta > 0$ such that $K_{\phi, 2} (\delta) \leq \eps$. 
Let $\mu$ be the measure $\ind_{\D \setminus \overline{D (0, 1 - \delta)}}\,. {\cal A}_\phi$. One has 
$K_{\mu, 2} (1) \leq 2 K_{\phi, 2} (\delta)$, because, for $\delta < h \leq 1$, the intersection of a window of size $h$ with the annulus 
$\{z \in \D\, ;\ 1 - \delta < |z| < 1\}$ can be covered by less than $2(h/\delta)$ windows of size $\delta$. \par 
Now, the Bergman version of Carleson's Theorem (see \cite{Stegenga}, proof of  
Theorem~1.2, bottom of the page 117, with $\alpha = 1/2$) says that:
\begin{displaymath}
\int_\D  |f (z)|^p \,d\mu \leq C_p\, K_{\mu, 2} (1)\, \| f \|_{{\mathfrak B}^p}^p 
\end{displaymath} 
for every $f \in {\mathfrak B}^p$. Hence: 
\begin{align*}
\int_\D |f_n \circ \phi|^p\,d{\cal A} 
& = \int_{\overline{\D (0, 1 - \delta)}} |f_n (z)|^p\,d{\cal A}_\phi (z) 
+ \int_{\D \setminus \overline{D (0, 1 - \delta)}} |f_n (z)|^p\,d{\cal A}_\phi (z) \\
& \leq \eps + 2 C_p \, \eps\,,
\end{align*}
for $n$ large enough, since $(f_n)$ converges uniformly to $0$ on $\overline{\D (0, 1 - \delta)}$. It follows that 
$\| C_\phi (f_n) \|_{{\mathfrak B}^p}$ converges to $0$.\qed
\goodbreak

\section{Hardy-Orlicz and Bergman-Orlicz spaces}\label{H-O and B-O spaces}

\subsection{Orlicz spaces}

An \emph{Orlicz function} is a function $\Psi \colon [0, \infty) \to [0, \infty)$ which is positive, non-decreas\-ing,  
convex and such that $\Psi (0) = 0$, $\Psi (x) > 0$ for $x > 0$ and 
$\Psi (x) \mathop{\longrightarrow}\limits_{x \to \infty} \infty$. 
\par\smallskip

\noindent{\bf Examples.} $\Psi (x) = x^p$; $\Psi (x) = x^p \log (x + 1)$, $1 \leq p < \infty$; 
$\Psi (x) = \e^{x^q} - 1$; $\Psi (x) = \exp [\big( \log (x + 1) \big)^q] - 1$, $q \geq 1$. 
\par\medskip

If $(S, {\cal T}, \mu)$ is a finite measure space, the \emph{Orlicz space} $L^\Psi (\mu)$ is the space of classes of 
measurable functions $f \colon S \to \C$ such that, for some $C > 0$:
\begin{displaymath}
\int_S \Psi (|f|/C) \,d\mu < \infty\,.
\end{displaymath}
The norm is defined by:
\begin{displaymath}
\| f \|_\Psi = \inf \{C \,;\ \int_S \Psi (|f|/C) \,d\mu \leq 1 \}\,.
\end{displaymath}
For $\Psi (x) = x^p$, we get the classical Lebesgue space: $L^\Psi (\mu) = L^p (\mu)$.

\subsection{Hardy-Orlicz and Bergman-Orlicz spaces}

We define the Bergman-Orlicz space ${\mathfrak B}^\Psi$ by:
\begin{displaymath}
{\mathfrak B}^\Psi = \{f \in L^\Psi (\D, {\cal A})\,;\ f \text{ analytic} \}\,,
\end{displaymath}
with the norm $\| f \|_{{\mathfrak B}^\Psi} = \| f \|_{L^\Psi (\D, {\cal A})}$.\par
The Hardy-Orlicz space $H^\Psi$ can be defined as in the classical case (see \cite{CompOrli}, Definition~3.2), but 
is is easier to define it by:
\begin{displaymath}
H^\Psi = \{f \in H^1\,;\ f^\ast \in L^\Psi (\T, m) \}\,,
\end{displaymath}
with the norm $\|f\|_{H^\Psi} = \| f^\ast \|_{L^\Psi (m)}$.\par
\bigskip

As in the classical case (because $\Psi (|f|)$ is subharmonic), Littlewood's subordination principle implies that 
every analytic $\phi \colon \D \to \D$ induces bounded composition operators 
$C_\phi \colon {\mathfrak B}^\Psi \to {\mathfrak B}^\Psi$ and $C_\phi \colon H^\Psi \to H^\Psi$.

\subsection{Compactness}

Theorem~\ref{compactness classical} has the following generalization (\cite{CompOrli}, Theorem~4.11) and 
\cite{LLQR-B}, Theorem~2.5):
\begin{theorem}\label{Carleson-Orlicz}
Let $\mu$ be a finite positive measure on $\D$, and assume that the identity maps 
$I_\mu \colon H^\Psi \to L^\Psi (\mu)$ and $J_\mu \colon {\mathfrak B}^\Psi \to L^\Psi (\mu)$ are defined. 
Then:\par\smallskip
\noindent $1)$ $\displaystyle \lim_{h \to 0} \frac{\Psi^{-1} (1/h)}{\Psi^{-1} [1 / h K_\mu (h)]} = 0 
\quad \Rightarrow \quad I_\mu \text{ compact} \quad \Rightarrow \quad 
\lim_{h \to 0} \frac{\Psi^{-1} (1/h)}{\Psi^{-1} [1/ \rho_\mu (h)]} = 0$;
\par\smallskip
\noindent $2)$ $\displaystyle \lim_{h \to 0} \frac{\Psi^{-1} (1/h^2)}{\Psi^{-1} [1 / h^2 K_{\mu, 2} (h)]} = 0 \!
\quad\! \Rightarrow \quad  J_\mu \text{ compact} \!\quad \Rightarrow \quad \! 
\lim_{h \to 0} \frac{\Psi^{-1} (1/h^2)}{\Psi^{-1} [1/ \rho_\mu (h)]} = 0$,
\par\smallskip
\noindent where $\displaystyle K_\mu (h) = \sup_{0 < t \leq h} \frac{\rho_\mu (t)}{t}$ and 
$\displaystyle K_{\mu,2} (h) = \sup_{0 < t \leq h} \frac{\rho_\mu (t)}{t^2}\,\cdot$
\end{theorem}

Actually, the sufficient conditions imply the existence of the identity maps $I_\mu$ and $J_\mu$.\par
\smallskip

The conditions in $1)$, resp. $2)$, are equivalent if $\Psi$ is ``regular'' (namely, if $\Psi$ satisfies the condition $\nabla_0$, whose definition is given after 
Proposition~\ref{equiv avec nabla} below), as $\Psi (x) = x^p$, for which case both  conditions in $1)$ read as 
$\displaystyle \frac{\rho_\mu (h)}{h} \mathop{\longrightarrow}_{h \to 0} 0$ and in $2)$ as 
$\displaystyle \frac{\rho_{\mu, 2} (h)}{h^2} \mathop{\longrightarrow}_{h \to 0} 0$, but examples show that there is 
no equivalence in general (see \cite{CompOrli}, pp. 50--54 and \cite{LLQR-B}, \S~2).\par

\subsection{Compactness for composition operators}

Nevertheless, for composition operators $C_\phi$, the following theorem (see \cite{CompOrli}, Theorem~4.19 and 
\cite{LLQR-B}, Theorem~3.1), which is one of the main result of this survey, says that:
\begin{displaymath}
K_\phi (h) \approx \rho_\phi (h)/h \quad \text{and} \quad 
K_{\phi, 2} (h) \approx \rho_{\phi, 2} (h)/ h^2.
\end{displaymath}

\begin{theorem}\label{homogene}
There is a universal constant $C > 0$ such that, for every analytic self-map $\phi\colon \D \to \D$, one has, for every 
$0 < \eps < 1$ and $h > 0$ small enough, for every $\xi \in \T$:\par\smallskip
$1)$ \qquad $m [\phi^\ast \in W (\xi, \eps h)] \leq C\,\eps\, m [\phi^\ast \in W (\xi, h)]$; \par\smallskip
$2)$ \qquad ${\cal A} [\phi\, \in W (\xi, \eps h)] \,\,\leq C\,\eps^2 \, {\cal A}[\phi \in W (\xi, h)]$. 
\end{theorem}

For fixed $h$, this expresses that the measure $m_\phi$ is a Carleson measure and ${\cal A}_\phi$ is a 
$2$-Carleson measure. The theorem says that this is true at ``all scales''.\par\smallskip
 
It follows that:

\begin{theorem}\label{compacite operateurs composition}
For every analytic self-map $\phi \colon \D \to \D$, one has:\par
\smallskip
$1)$ \quad $C_\phi \colon H^\Psi \to H^\Psi$ compact $\Longleftrightarrow$ \quad 
$\displaystyle \lim_{h \to 0} \frac{\Psi^{-1} (1/h)}{\Psi^{-1} [1/\rho_\phi (h)]} = 0$;\par 
\smallskip
$2)$ \quad $C_\phi \colon {\mathfrak B}^\Psi \to {\mathfrak B}^\Psi$ compact $\Longleftrightarrow$ \quad 
$\displaystyle \lim_{h \to 0} \frac{\Psi^{-1} (1/h^2)}{\Psi^{-1} [1/\rho_{\phi, 2} (h)]} = 0$. 
\end{theorem}

Let us give a very vague idea of the proof of Theorem~\ref{homogene}, in the Bergman case (the proof in the Hardy case, 
though following the same ideas, is different, and actually more difficult). By setting $f = h/ (1 - \phi)$, it 
suffices to show that
\begin{equation}\label{un}
\qquad \quad {\cal A} (\{ |f| > \lambda\}) \leq \frac{K}{\lambda^2}\,{\cal A} (\{|f| > 1\}) \qquad \text{for } 
\lambda \geq \lambda_0 > 0\,,
\end{equation}
for every analytic function $f \colon \D \to \Pi^+ = \{z\in \C\,;\ \Re z > 0\}$ such that $|f (0)| \leq \alpha_0$, for 
some $\alpha_0 > 0$. But the fact that ${\cal A}_\phi$ is a $2$-Carleson measure writes:
\begin{equation}\label{deux}
{\cal A} (\{ |f| > \lambda\}) \leq \frac{C}{\lambda^2}\,|f (0)|\,, \qquad \lambda > 0\,,
\end{equation}
where $f = h/ (1 - \phi)$. One has hence to replace $|f (0)|$ in the majorization by ${\cal A} (\{|f| > 1\})$. 
To that effect, one splits the disk $\D$ into pieces which are ``uniformly conform'' to $\D$ and on which we can 
use \eqref{deux}. However, it is far from being so easy, and we refer to \cite{LLQR-B} (and \cite{CompOrli} for the 
Hardy case) for the details.
\goodbreak

\section{Compactness on $H^\Psi$ versus compactness on ${\mathfrak B}^\Psi$}\label{section comparaison}

Thanks to Theorem~\ref{compacite operateurs composition}, in order to compare the compactness of the composition operator 
$C_\phi$ on $H^\Psi$ and on ${\mathfrak B}^\Psi$, we have to compare $\rho_\phi (h)$ and $\rho_{\phi, 2} (h)$. But 
if one reads their definitions:
\begin{align*}
\rho_\phi (h) & = \sup_{|\xi| = 1} m [\phi^\ast \in W (\xi, h)] \\
\noalign{\noindent \text{and}}\\
\rho_{\phi, 2} (h) & = \sup_{|\xi| = 1} {\cal A} [ \phi \in W (\xi, h)]\,,
\end{align*}
that does not seem straightforward. We shall compare them in an indirect way, by using the Nevanlinna counting 
function.

\subsection{Nevanlinna counting function}

The \emph{Nevanlinna counting function} counts how many pre-images each element has, with a weight which decreases 
when this pre-image approaches $\partial \D$. Namely:
\begin{displaymath}
N_\phi (w) = \sum_{\phi (z) = w} \log \frac{1}{|z|} 
\end{displaymath}
for $w \in \phi (\D)$ and $w \neq \phi (0)$. One sets $N_\phi (w) = 0$ for the other $w \in \D$.\par\smallskip

Our second main theorem asserts that the Nevanlinna counting function of $\phi$ is equivalent to its Carleson 
function (see \cite{LLQR-N}). 
\goodbreak 

\begin{theorem}
There exists a universal constant $C > 1$ such that, for every analytic self-map $\phi\colon \D \to \D$ and for 
$h > 0$ small enough, we have:
\begin{displaymath}
\frac{1}{C} \sup_{w \in W (\xi, h/C)} N_\phi (w) \leq m_\phi [W (\xi, h)] \leq 
\frac{C}{{\cal A}[W (\xi, Ch)]} \int_{W (\xi, Ch)} \hskip -15 pt N_\phi (z) \,d{\cal A} (z)\,.
\end{displaymath}
\end{theorem}

Now, if one defines the Nevanlinna function of order $2$ by:
\begin{displaymath}
N_{\phi, 2} (w) = \sum_{\phi (z)= w} \bigg[\log \frac{1}{|z|}\bigg]^2 
\end{displaymath}
for $w \in \phi (\D)$ and $w \neq \phi (0)$ and $N_{\phi, 2} (w) = 0$ for the other $w \in \D$, one has easily 
(see \cite{Shapiro-Annals}, Proposition~6.6):
\begin{displaymath}
N_{\phi, 2} (w) = 2 \int_0^1 N_\phi (r, w)\, \frac{dr}{r}\,,
\end{displaymath}
where $\displaystyle N_\phi (r, w) = \sum_{\phi (z) = w, |z| < r} \log \frac{r}{|z|}$ is the restricted Nevanlinna 
function. But, since $N_\phi (r, w) = N_{\phi_r} (w)$ with $\phi_r (z) = \phi (rz)$, one gets, by integrating in 
polar coordinates:
\begin{corollary}\label{equiv 2-Carleson et 2-Nevanlinna}
There is some universal constant $C > 1$ such that, for every analytic self-map $\phi\colon \D \to \D$ and for 
$h > 0$ small enough:
\begin{displaymath}
\frac{1}{C}\, \rho_{\phi, 2} (h/C) \leq \sup_{|w| \geq 1 - h} N_{\phi, 2} (w) \leq C\, \rho_{\phi, 2} (C h)
\end{displaymath}
\end{corollary} 

But now, it is easy to compare $N_\phi$ and $N_{\phi, 2}$: one has: 
\begin{displaymath}
N_{\phi, 2} (w) \leq [N_\phi (w)]^2 
\end{displaymath}
(simply because the $\ell_2$-norm is smaller than the $\ell_1$-norm). Therefore:
\begin{theorem}
There is a universal constant $C > 1$ such that, for every analytic self-map $\phi\colon \D \to \D$, one has for 
$h > 0$ small enough and every $\xi\in \T$:
\begin{displaymath}
{\cal A_\phi} [W (\xi, h)] \leq C\, \big( m_\phi [W (\xi, Ch)] \big)^2
\end{displaymath}
\end{theorem}

We can now compare the compactness on $H^\Psi$ and on ${\mathfrak B}^\Psi$.

\begin{theorem}
Under some condition on $\Psi$, one has:
\begin{displaymath}
C_\phi \colon H^\Psi \to H^\Psi \quad \text{compact} \quad \Longrightarrow \quad 
C_\phi \colon {\mathfrak B}^\Psi \to {\mathfrak B}^\Psi \quad \text{compact}.
\end{displaymath}
\end{theorem}

This condition has not a very nice statement:
\begin{displaymath}
\forall A > 0, \ \exists x_A > 0, \ \exists B \geq A: \qquad 
\Psi [A \Psi^{-1} (x^2)] \leq \big(\Psi [B \Psi^{-1} (x)]\big)^2, \qquad x \geq x_A 
\end{displaymath}
(though, setting $\chi_A (x) = \Psi [A \Psi^{-1} (x)]$, it writes better $\chi_A (x^2) \leq [\chi_B (x)]^2$), 
but it is satisfied in many cases:
\begin{itemize}
\item[-] if $\Psi$ grows moderately; namely, satisfies the condition $\Delta_2$, {\it i.e.} 
$\Psi (2x) \leq C\, \Psi (x)$ for $x$ large enough; this is the case for $\Psi (x) = x^p$, and we recover the 
classical case of Corollary~\ref{comparaison cas classique};
\item[-] if $\Psi$ grows quickly; namely, satisfies the condition $\Delta^2$, {\it i.e.} for some $\alpha > 1$, 
$[\Psi (x)]^2 \leq \Psi (\alpha x)$ for $x$ large enough; for instance if $\Psi (x) = \e^{x^q} - 1$, $q \geq 1$; 
\item[-] but also for $\Psi(x) = \exp [ \big( \log (x + 1) \big)^2] - 1$.
\end{itemize}

Nevertheless (\cite{LLQR-B}, Theorem~4.2): 

\begin{theorem}
There exists a symbol $\phi$ and an Orlicz function $\Psi$ such that the composition operator 
$C_\phi \colon H^\Psi \to H^\Psi$ is compact, and moreover in all Schatten classes $S_p  (H^2)$, $p > 0$, whereas 
$C_\phi \colon {\mathfrak B}^\Psi \to {\mathfrak B}^\Psi$ is \emph{not} compact.
\end{theorem}

The symbol $\phi$ is a conformal map from $\D$ onto the domain $G$, represented on the picture, delimited by three 
circular arcs of radii $1/2$.

\begin{figure}[t]
\begin{center}
\includegraphics[width=0.3\linewidth]{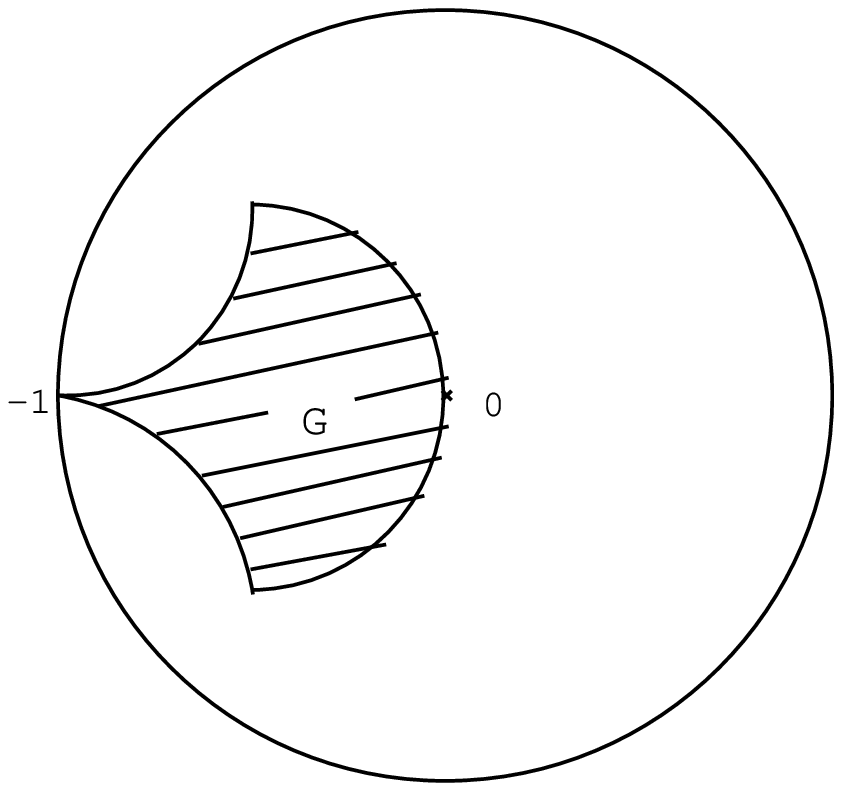}
\end{center}
\end{figure}

\medskip\goodbreak
The Carleson function of $\phi$ is ``small'' whereas its Carleson function of order $2$ is ``small'':
\begin{align*}
\rho_\phi (h) & \leq C\, \e^{- \pi/4 h} \\
\rho_{\phi, 2} (h) & \geq (1/C)\, \e^{- \pi /h}
\end{align*}

Now, it remains to construct a concave and piecewise linear function $F$ (and $\Psi$ will be $F^{-1}$) in such way 
that:
\begin{displaymath}
\lim_{x \to \infty} \frac{F (x)}{F (\e^{\pi x/4})} = 0 \quad \text{and} \quad 
\limsup_{x \to \infty} \frac{F (x^2)}{F (\e^{\pi x})} > 0 \,,
\end{displaymath}
in order that $\displaystyle \lim_{h \to 0} \frac{\Psi^{-1} (1/h)}{\Psi^{-1} [1/\rho_\phi (h)]} = 0$, but 
$\displaystyle \limsup_{h \to 0} \frac{\Psi^{-1} (1/h^2)}{\Psi^{-1} [1/\rho_{\phi, 2} (h)]} > 0$. \qed
\par\bigskip

\noindent{\bf Remark.} Let us stress that the compactness of $C_\phi \colon H^\Psi \to H^\Psi$ always implies the 
compactness of $C_\phi \colon {\mathfrak B}^{\Psi^2} \to {\mathfrak B}^{\Psi^2}$. Indeed, if 
$\tilde \Psi (x) = [\Psi (x)]^2$, then $\tilde \Psi^{-1} (t^2) = \Psi^{-1} (t)$, so we get 
$\tilde \Psi^{-1} (1/h^2) / \tilde \Psi^{-1} \big(1/ \nu_{\phi, 2} (h) \big) \leq 
\Psi^{-1} (1/h) / \Psi^{-1} \big(1/ \nu_\phi (h) \big)$, since one has $\nu_{\phi, 2} (h) \leq [\nu_\phi (h)]^2$, 
where $\nu_\phi (h) = \sup_{|w| \geq 1 - h} N_\phi (w)$.\par

\section{Final remarks}

\subsection{Modulus of the symbol and compactness on Bergman-Orlicz spaces}

What about the compactness of $C_\phi$ on the Bergman-Orlicz spaces and the behaviour of the 
modulus of its symbol $\phi$? \par
We prove in \cite{CompOrli}, Theorem~5.7, that the compactness of 
$C_\phi \colon {\mathfrak B}^\Psi \to {\mathfrak B}^\Psi$ implies that:
\begin{equation}\label{condition module Bergman}
\lim_{|a| \to 1} 
\frac{\Psi^{-1} \bigg[\displaystyle \frac{1}{(1 - |\phi (a)|)^2}\bigg]}
{\Psi^{-1} \bigg[\displaystyle \frac{1}{(1 - |a|)^2}\bigg]} 
= 0\,.
\end{equation}
The proof is essentially the same as the proof of Proposition~\ref{Shapiro-Taylor} given in 
Section~\ref{some comments}. But it follows also from Theorem~\ref{compacite operateurs composition} and 
Corollary~\ref{equiv 2-Carleson et 2-Nevanlinna}. Indeed, these results have the following consequence:
\begin{theorem}
For every analytic self-map $\phi \colon \D \to \D$ and for every Orlicz function $\Psi$, the composition 
operator $C_\phi \colon {\mathfrak B}^\Psi \to {\mathfrak B}^\Psi$ is compact if and only if:
\begin{equation}\label{petit o Nevanlinna}
\lim_{h \to 0} \frac {\Psi^{-1} (1/h^2)} {\Psi^{-1} \big(1/ \nu_{\phi, 2} (h) \big) } = 0\,,
\end{equation}
where $\nu_{\phi, 2} (h) = \sup_{|w| \geq 1 - h} N_{\phi, 2} (w) $.
\end{theorem}
Now, since $N_{\phi, 2} \big( \phi (z) \big) \geq \big( \log (1 /|z|) \big)^2 \geq ( 1 - |z|)^2$, we get that the 
compactness of $C_\phi \colon {\mathfrak B}^\Psi \to {\mathfrak B}^\Psi$ implies  
\eqref{condition module Bergman}. But Theorem~\ref{petit o Nevanlinna} gives a partial converse:
\begin{theorem}\label{univalent}
For every \emph{univalent} (or more generally, boundedly-valent) analytic self-map $\phi \colon \D \to \D$ and for 
every Orlicz function $\Psi$, the composition operator $C_\phi \colon {\mathfrak B}^\Psi \to {\mathfrak B}^\Psi$ is 
compact if and only if one has \eqref{condition module Bergman}.
\end{theorem}

Let us recall that $\phi$ bounded-valent means that there is an integer $L \geq 1$ such that the equation 
$\phi (z) = w$ has at most $L$ solution(s) in $\D$, for every $w \in \D$; we then have 
$N_{\phi, 2} (w) \leq L (1 - |z|)$, where $\phi (z) = w$, with $|z| > 0$ minimal, and Theorem~\ref{univalent} 
follows.\qed
\par\medskip 

Compactness of $C_\phi$ also is equivalent to \eqref{condition module Bergman}, for every symbol $\phi$, if one 
adds a condition on the Orlicz function $\Psi$.

\begin{proposition}\label{equiv avec nabla}
Assume that the Orlicz function $\Psi$ satisfies the condition $\nabla_0$. Then the composition operator 
$C_\phi \colon {\mathfrak B}^\Psi \to {\mathfrak B}^\Psi$ is compact if and only if \eqref{condition module Bergman} 
holds.
\end{proposition}

Condition $\nabla_0$ is the ``regularity'' condition mentioned after Theorem~\ref{Carleson-Orlicz}, which gives 
the equivalence between the necessary and sufficient conditions. $\Psi$ satisfies $\nabla_0$ if 
(see \cite{CompOrli}, Definition~4.5 and Proposition 4.6): there is $x_0 > 0$ such that, for every $B > 1$, there exists $c_B > 1$ such that:
\begin{equation}\label{nabla}
\frac{\Psi (B x)}{\Psi (x)} \leq \frac{\Psi (c_B B y)}{\Psi (y)} \qquad \text{for } x_0 \leq x \leq y.
\end{equation} 
Let us point out that it is satisfied, in particular, when $\log \Psi (\e^x)$ is convex, or if $\Psi$ satisfies the 
condition $\Delta^2$ (\cite{CompOrli}, Proposition~4.7).
\par\medskip

\noindent{\bf Proof.} We may assume that $\phi (0) = 0$. Then Littlewood's inequality writes 
$N_\phi (w) \leq \log 1/|w|$, $w \neq 0$, and gives, for some constant $C > 1$:
\begin{displaymath}
N_{\phi, 2} (w) \leq C \sup_{\phi (z) = w} (1 - |z|) \sum_{\phi (z) = w} \log \frac{1}{|z|} 
\leq C^2(1 - |w|)\,\sup_{\phi (z) = w} (1 - |z|)\,.
\end{displaymath}

By hypothesis, for every $A >0$, we have, with $w = \phi (z)$:
\begin{displaymath}
\Psi^{-1}  \bigg[ \frac{1}{( 1 - |z|)^2} \bigg] \geq A\, \Psi^{-1} 
\Bigg[\displaystyle \frac{1}{( 1 - |\phi (z)|)^2} \Bigg]
\end{displaymath}
for $|z|$ close enough to $1$. With $h = 1 - |w|$, this writes:
\begin{displaymath}
1 - |z| \leq 1/ \sqrt{\Psi \big[A \Psi^{-1} (1/h^2)\big]}\,.
\end{displaymath}
We get hence:
\begin{displaymath}
\nu_{\phi,2} (h) \leq \frac{C^2 h}{\sqrt{\Psi \big[A \Psi^{-1} (1/h^2)\big]}} 
\end{displaymath}
if $h > 0$ is small enough. It follows that we shall have:
\begin{displaymath}
\lim_{h \to 0} \frac{\Psi^{-1} (1/h^2)}{\Psi^{-1} \big( 1/\nu_{\phi,2} (h) \big) } = 0
\end{displaymath}
if for every $B > 1$, we can find $A > 0$ such that:
\begin{displaymath}
\Psi^{-1} \bigg[ \frac{ \sqrt{\Psi [ A \Psi^{-1} (1/h^2) ]}} {C^2 h} \bigg] \geq B \Psi^{-1} (1/h^2)\,.
\end{displaymath}
Setting $x = \Psi^{-1} (1/h^2)$, it suffices to have:
\begin{equation}\label{condition Psi}
\big[\Psi (B x)\big]^2 \leq \Psi (x)\, \Psi \big((A/C^4)x \big)\,,
\end{equation}
for $x > 0$ big enough, since $\Psi \big( (A/C^4) x \big) \leq \Psi (Ax)/ C^4$ by convexity.\par\smallskip

When $\Psi \in \nabla_0$, \eqref{nabla} gives \eqref{condition Psi} with $y = Bx$ and $A = C^4 c_B B^2$.\qed

\subsection{Blaschke products}

We may ask about the converse implication: does the compactness of $C_\phi$ on 
${\mathfrak B}^\Psi$ imply its compactness on $H^\Psi$? Even in the Hilbertian case ${\mathfrak B}^2-H^2$, this is 
not the case (see \cite{Shapiro}, pp. 183--185): there is a Blaschke product $B$ (whose associated composition 
operator is an isometry from $H^2$ into itself) with no angular derivative, so $C_B$ is compact on ${\mathfrak B}^2$. Another 
example (a Blaschke product also) is given in \cite{CompOrli}, when $\Psi (x) = \e^{x^2} - 1$. More generally:

\begin{theorem}
For every Orlicz function $\Psi$ satisfying the condition $\nabla_0$, there is a Blaschke product $B$, whose associated composition operator  
$C_B$ is an isometry from $H^\Psi$ into itself, but such that $C_B$ is compact on ${\mathfrak B}^\Psi$. 
\end{theorem}
\noindent{\bf Proof.} Indeed (\cite{LLQR-Rev}, Theorem~3.1), for every function 
$\delta \colon (0, 1) \to (0, 1/2]$ such that $\delta (t) \mathop{\longrightarrow}\limits_{t \to 0} 0$, 
there is a Blaschke product $B$ such that:
\begin{equation}\label{Blaschke}
\qquad \qquad 1 - |B(z)| \geq \delta (1 - |z|), \qquad \text{for all } z\in \D\,.
\end{equation}
By replacing $B (z)$ by $z B (z)$, we may assume that $B (0) = 0$ (note that $1 - |z B (z)| \geq 1 - |B (z)|$). Then $C_B$ is an isometry from 
$H^\Psi$ into itself. Indeed, one can see (\cite{Nordgren}, Theorem~1) that if $\phi$ is an inner function, then the pull-back measure $m_\phi$ is 
equal to $P_a.m$, when $P_a$ is the Poisson kernel at $a = \phi (0)$. When $\phi (0) = 0$, one has $m_\phi = m$ and $C_\phi$ is an isometry. \par
Hence, taking, for $t > 0$ small enough:
\begin{displaymath}
\delta (t) = \frac{1} {\sqrt{\Psi \big[\sqrt{ \Psi^{-1} (1/t^2)} \big]}}\,,
\end{displaymath}
we get 
\begin{displaymath}
\frac{\Psi^{-1} \bigg[\displaystyle \frac{1}{(1 - |B (z)|)^2}\bigg]}
{\Psi^{-1} \bigg[\displaystyle \frac{1}{(1 - |z|)^2}\bigg]} 
\leq \frac{1}{\sqrt{\Psi^{-1} \bigg[ \displaystyle \frac{1}{(1 - |z|)^2} \bigg]}} \, \raise 1 pt \hbox{,}
\end{displaymath}
and we get \eqref{condition module Bergman}. Hence, assuming that $\Psi$ satisfies the condition $\nabla_0$, 
$C_B$ is compact on ${\mathfrak B}^\Psi$, by Proposition~\ref{equiv avec nabla}.\qed 
\par\bigskip

\noindent{\bf Remark.} Another, but different, survey on this topic can be found in \cite{Herve}


\nobreak
\vbox{\small \noindent{\it
{\rm Daniel Li}, Univ Lille Nord de France F-59\kern 1mm 000 LILLE, 
U-Artois, Laboratoire de Math\'ematiques de Lens EA~2462, 
F\'ed\'eration CNRS Nord-Pas-de-Calais FR~2956, 
F-62\kern 1mm 300 LENS,
Facult\'e des Sciences Jean Perrin,
Rue Jean Souvraz, S.P.\kern 1mm 18, FRANCE \\ 
daniel.li@euler.univ-artois.fr }}

\end{document}